\documentclass[12pt]{article}
\usepackage{amsmath}
\usepackage{amsfonts}
\usepackage{amssymb}
\usepackage{graphicx}
\usepackage{epsfig}
\font\of=msbm10 scaled 1200
\def\R{\mbox{\of R}}
\def\C{\mbox{\of C}}
\def\Z{\mbox{\of Z}}
\def\N{\mbox{\of N}}
\def\P{\mbox{\of P}}
\title{ Perturbations of  quadratic centers of genus one}
 \author{S\'{e}bastien Gautier, Lubomir Gavrilov \\
 \normalsize \it Institut de Math\'{e}matiques,
 \normalsize \it Universit\'{e} de Toulouse\\
   \normalsize \it 31062 Toulouse cedex 9, \normalsize \it   France
   \\ Iliya D. Iliev\\
\normalsize \it Institute of Mathematics, Bulgarian Academy of Sciences\\
\normalsize \it Bl. 8, 1113 Sofia, Bulgaria }

\begin{document}

\maketitle
 \hoffset =-1truecm \voffset =-2truecm
 \font\of=msbm10 scaled 1200

\begin{abstract}
We propose a program for finding the cyclicity of period annuli of quadratic
systems with centers of genus one.  As a first step, we classify all such
systems and determine the essential one-parameter quadratic perturbations
which produce the maximal number of limit cycles.  We compute the associated
Poincar\'{e}-Pontryagin-Melnikov functions whose zeros control the number of
limit cycles.  To illustrate our approach, we determine the cyclicity of the
annuli of two particular reversible systems.  \end{abstract}

\vspace{2ex}\noindent {\large\bf 0. Introduction}

\vspace{2ex} \noindent
As well known, there are four types of planar quadratic systems with a center:
1) Hamiltonian, 2) reversible, 3) generalized Lotka-Volterra, and 4) of
codimension four.  The first integral in the Hamiltonian case is a cubic
polynomial in $(x,y)$ whose generic level sets are elliptic curves. It is easy
to observe that the generic level sets of the first integral of the codimension
four
case are elliptic curves too (see the end of Section 4). Our main efforts will
be devoted to the remaining two cases. It has been recently proved by one of
the authors (S.G.) that there are 18 classes of reversible centers, 6 classes
of reversible Lotka-Volterra centers and 5 classes of generic (nonreversible)
Lotka-Volterra centers whose phase portraits contain only elliptic curves
(possibly, a finite number of them reducible) \cite{3}.  They are given by
codimension-one or codimension-two algebraic sets in the space of all centers
from the corresponding type, see Theorems 1 and 2 below. Throughout the paper,
by "genus" we mean the genus of the compactified and normalized generic phase
curves. An algebraic phase curve is generic if it does not contain a singular
point of the vector field in its closure. The centers whose (generic complexified)
periodic orbits are elliptic curves will be called \emph{centers of genus one}.
We note that even a quadratic system can have a center of arbitrarily big genus.

Once we know that the first integral of a given planar system with a center
defines elliptic curves, we could raise the next question: how many limit
cycles could be produced in the phase portrait under small quadratic (or even
polynomial) perturbations?  \emph{The purpose of the paper is to present a
program for solving this problem.} In what follows we restrict our attention to
the limit cycles which bifurcate from open period annuli (we do not consider
graphics).  It turns out that, instead of multi-parameter perturbations, it is
enough to consider suitable \emph{one-parameter} small quadratic perturbations,
see \cite{gav1}.  The one-parameter perturbations which produce the maximal
number of limit cycles in the quadratic case, called \emph{essential
perturbations}, together with the corresponding generating functions of limit
cycles (or also Poincar\'{e}-Pontryagin-Melnikov functions) were determined in
\cite{5}.

The first part of our program is to adapt the results of \cite{5} to our case.
The result is a complete list of such essential perturbations of quadratic
systems with centers of genus one, together with the corresponding generating
functions.  The list is presented in Section 3 (in the reversible case) and in
Section 4 (in the Lotka-Volterra case).

The second part of the program is to study the zeros of the generating
functions $I(t)$ found in Sections 3 and 4.  The fact that each
$I(t)=\int_{\{H=t\}}\omega$ is a complete elliptic integral from a rational
one-form $\omega$, over the level sets $\{H=t\}\subset\C^2$ which are elliptic
curves, allows one to apply all related facts from algebraic geometry in order
to estimate the number of zeros of $I(t)$ and thus to set up some upper (or
lower) bounds on the number of bifurcating limit cycles in the system, see e.g.
\cite{8,11,hor,4,gav4}. This part of our program is illustrated in Section 5,
where we study the cyclicity of the period annuli for two types of quadratic
reversible systems with centers of genus one.

The paper is organized as follows.  In Section 1 we determine all reversible
centers with phase portrait formed by elliptic curves. In Section 2 we discuss
the same question for the Lotka-Volterra centers (reversible or not).  These
results were previously proved in \cite{3}.  For convenience of the reader we
present here almost self-contained proofs adapted for the purposes
of the present paper. In Sections 3 and 4 we determine for each of the cases
the corresponding generating function, the complete elliptic integral $I(t)$
which is the leading term in the expansion of the first return mapping,
respectively for the reversible and Lotka-Volterra cases. We also present
there (as Conjectures 1 and 2) the expected exact upper bounds for the number
of zeros of all generating functions. In Section 5 we
use the geometric properties of the elliptic fibration determined by the first
integral of the quadratic system, in order to determine, for two of the cases,
the number of the zeros of $I(t)$.  From this we deduce an exact result: the
cyclicity of the annuli under consideration is two (Theorem 3).

\newpage
\vspace{2ex} \noindent {\large\bf 1. The reversible case}

\vspace{2ex} \noindent The general first integral of a reversible system
with a center at the origin $$\dot{z}=-iz+az^2+2|z|^2+b\bar{z}^2,\quad
a,b\in \R,\quad z=x+iy\eqno(1.1)$$ is given by $$\textstyle
H(x,y)=X^\lambda(\frac12 y^2+AX^2+BX+C)\eqno(1.2)$$ where $X=1+2(a-b)x$ and
$\lambda,A,B,C$ are explicit rational functions of the parameters $a,b$
(see formula (1.7) below).  Moreover, one has $\lambda\neq 0, -1, -2.$
Along with the elliptic Hamiltonian $H=\frac12y^2+\frac12x^2-\frac43x^3$
(corresponding to $a=b=-1$ in (1.1)),
formula (1.2) contains all the cases for which $\{H=t\}$ is (possibly) an
algebraic curve of genus 1.  Below, we are going to determine all these
cases.

We begin with the observation that it is enough to consider only the case
when $\lambda< -1.$ Indeed, if $\lambda> -1$, by the bi-rational change
$(X,y)=(1/X_1,Y_1/X_1)$ the function (1) reduces to
$H_1=X_1^{-2-\lambda}(\frac12 Y_1^2+CX_1^2+BX_1+A)$ with $-2-\lambda<-1.$
For this reason below we will consider the case when
$$\lambda=-\frac{p}{q}\quad \mbox{\rm with}\quad p,q
\in\N,\;\;p>q,\;\;p\neq 2q,\;\;gcd(p,q)=1$$
where $gcd(p,q)$ denotes the greatest common divisor of $p$ and $q$.

\vspace{2ex} \noindent
1.  Assume first that $A\neq 0$, $C\neq 0$.  Then $\{H=t\}$ after taking
$X=x^q$ (this is an isomorphic map \cite{3}, Lemma 1)
 reads $$\textstyle \frac12 y^2=-Ax^{2q}-Bx^q-C+tx^p.\eqno(1.3)$$ For
arbitrary $t$, the curve (1.3) has a genus 1 if and only if the polynomial on
the right hand side has a degree 3 or 4.  That is, either $p=3$ and $q=1,2$
or $p=4$ and $q=1.$

Consider now the degenerate cases when either $A=0$ or $C=0$ (let us note
that neither two of the coefficients $A,B,C$ can vanish simultaneously).

\vspace{2ex} \noindent
2.  Assume that $A=0.$ Then (1.3) is of genus 1 if and only if either $p=3,$
$q=1,2$ or $p=4,$ $q=1,3.$ Thus the unique new case compared to the case
$A\neq 0$ is $(p,q)=(4,3)$.

\vspace{2ex} \noindent
3.  Let $C=0.$ If $2q>p>q,$ by a bi-rational transformation
$(x,y)=(1/x_1,y_1/x_1^q)$ one can transform (1.3) into $$\textstyle \frac12
y^2=-A-Bx^q+tx^{2q-p}\eqno(1.4)$$ (here and below, we shall omit the
subscript 1).  As $q>2q-p$, the genus is 1 if and only if either $q=3$,
$p=4,5$ or $q=4$, $p=5,7$.  If $p>2q$ and $q$ is even, then the rational
change $y=y_1x^{q/2}$ in (1.3) yields $$\textstyle \frac12
y^2=-Ax^q-B+tx^{p-q}.\eqno(1.5)$$ As $p-q>q$, (1.5) is of genus 1 if and only
if $(p,q)=(5,2)$.  Finally, if $p>2q$ and $q$ is odd, the rational change
$y=y_1x^{(q-1)/2}$ in (1.3) yields $$\textstyle \frac12
y^2=-Ax^{q+1}-Bx+tx^{p-q+1}.\eqno(1.6)$$ As above, (1.6) is of genus one if and
only if $p-q+1$ equals 3 or 4, which is possible when either $(p,q)=(3,1)$
or $(p,q)=(4,1)$, cases that already have been obtained when we assumed
that $A\neq 0$.

\vspace{2ex} \noindent Thus we have obtained the complete list of
cases for which $\{H=t\}$, $H$ given by (1.2) and $t$ arbitrary, is
a curve of genus one (the right column contains the cases with
$\lambda>-1$):
$$\begin{array}{ll} {\rm (0)}\; H=\frac12y^2+Ax^2+Bx^3
&\mbox{\rm (the standard elliptic case)}\\ {\rm (i)}\; H=X^{-3}(\frac12
y^2+AX^2+BX+C) & {\rm (ii)}\; H=X(\frac12 y^2+CX^2+BX+A)\\
{\rm (iii)}\; H=X^{-3/2}(\frac12 y^2+AX^2+BX+C) & {\rm (iv)}\;
H=X^{-1/2}(\frac12 y^2+CX^2+BX+A)\\ {\rm (v)}\; H=X^{-4}(\frac12
y^2+AX^2+BX+C) & {\rm (vi)}\; H=X^2(\frac12 y^2+CX^2+BX+A)\\
{\rm (vii)}\; H=X^{-4/3}(\frac12 y^2+BX+C) & {\rm (viii)}\;
H=X^{-2/3}(\frac12 y^2+CX^2+BX)\\ {\rm (ix)}\;
H=X^{-4/3}(\frac12 y^2+AX^2+BX) & {\rm (x)}\; H=X^{-2/3}(\frac12
y^2+BX+A)\\ {\rm (xi)}\; H=X^{-5/3}(\frac12 y^2+AX^2+BX) &
{\rm (xii)}\; H=X^{-1/3}(\frac12 y^2+BX+A)\\ {\rm (xiii)}\;
H=X^{-5/4}(\frac12 y^2+AX^2+BX) & {\rm (xiv)}\;
H=X^{-3/4}(\frac12 y^2+BX+A)\\ {\rm (xv)}\; H=X^{-7/4}(\frac12
y^2+AX^2+BX) & {\rm (xvi)}\; H=X^{-1/4}(\frac12 y^2+BX+A)\\
{\rm (xvii)}\; H=X^{-5/2}(\frac12 y^2+AX^2+BX) & {\rm (xviii)}\;
H=X^{1/2}(\frac12 y^2+BX+A).
\end{array}$$
We should mention that in cases (iii) and (iv) above
it was assumed that $C\neq 0$.  If $C=0$, the curve has a genus
zero, see below.

 \vspace{1ex}
\noindent Let us now recall the exact formula of (1.2) from \cite{5}:
$$H(X,y)=X^{-\frac{a+b+2}{a-b}}\left(\frac{y^2}{2}+\frac{1}{8(a-b)^2}
\left(\frac{a+b-2}{a-3b-2}X^2
+2\frac{b-1}{b+1}X+\frac{a-3b+2}{a+b+2}\right)\right).\eqno(1.7) $$
This formula holds outside the lines $a=b$, $a+b+2=0$, $b=-1$, $a-3b-2=0$
(note the last three cases correspond to $\lambda=0,-1,-2$ in (1.2)). Except
for the three points $(a,b)=(-1,-1), (\pm2,0)$, on these lines the first
integrals contain exponents \cite{5} and hence their level sets are not algebraic
curves. For $(a,b)=(\pm2,0)$ the level sets are conic ovals. The lines
$a+b+2=-\lambda(a-b)$ in the $(a,b)$-plane, $\lambda\in\R$, together with
$a=b$, form the bundle of straight lines through the point $(-1,-1)$ which
corresponds to the standard elliptic case. Therefore using the above
results, we obtain:

\vspace{2ex} \noindent {\bf Theorem 1.} {\it The phase curves of
$(1.1)$ are algebraic curves of genus one if and only if one of the
conditions holds:} $$\begin{array}{ll} {\rm (r1)}\; a=2b+1 & {\rm
(r2)}\; a=-1 \;\mbox{\rm (the reversible Hamiltonian case)} \\
{\rm (r3)}\; a=5b+4 \;\;(b\neq-3) & {\rm (r4)}\; a=-3b-4
\;\;(b\neq-3) \\   {\rm (r5)}\; a=\frac53b+\frac23 & {\rm (r6)}\;
a=\frac13b-\frac23 \\ {\rm (r7)}\; (a,b)=(\frac52,-\frac12) & {\rm
(r8)}\; (a,b)=(-\frac72,-\frac12) \\ {\rm (r9)}\; (a,b)=(-8,-2) &
{\rm (r10)}\; (a,b)=(4,-2) \\ {\rm (r11)}\; (a,b)=(-17,-5) & {\rm
(r12)}\; (a,b)=(7,-5)
\\ {\rm (r13)}\; (a,b)=(-7,-\frac53) & {\rm (r14)}\;
(a,b)=(\frac{11}{3},-\frac53) \\ {\rm (r15)}\; (a,b)=(-23,-7) & {\rm
(r16)}\; (a,b)=(9,-7) \\ {\rm (r17)}\; (a,b)=(13,5) & {\rm (r18)}\;
(a,b)=(-3,5).  \end{array}$$
For completeness, below we add the list of all cases in (1.1) for which the
ovals $\{H=t\}$ are conic curves (ellipses).  These are
$$\begin{array}{l} {\rm (r19)}\; H=X^{-3/2}(\frac12 y^2+AX^2+BX) \;\;
\mbox{\rm when} \;\; (a,b)=(-11,-3)\\ {\rm (r20)}\; H=X^{-1/2}(\frac12
y^2+BX+A) \;\; \mbox{\rm when} \;\; (a,b)=(5,-3)\\ {\rm (r21)}\;
H=X^{-2}(\frac12 y^2+BX+C) \;\; \mbox{\rm when} \;\; (a,b)=(2,0)\\ {\rm
(r22)}\; H=\frac12 y^2+CX^2+BX \;\; \mbox{\rm when} \;\; (a,b)=(-2,0).
\end{array}$$
 Cases (r21) and (r22) (not included in (1.2) and (1.7)) are taken
from the full list of first integrals of (1.1), see e.g.  \cite{5}.
 The cases
(r19) and (r20) are obtained from (1.2) in the same way as above.
Note that (r20) and (r21) are the isochronous centers ${\cal S}_3$ and
${\cal S}_2$, respectively. By the way, the isochronous center
${\cal S}_4$ corresponds to $b=1$ in (r5). Perturbations of the
quadratic isochronous centers have been studied in \cite{12}.

\vspace{2ex}
 \noindent
 {\large\bf 2. The Lotka-Volterra case}

\vspace{2ex} \noindent A (generalized) Lotka-Volterra system with
a center at the origin has in complex coordinates the form
$$\dot{z}=-iz+Az^2+B\bar{z}^2,\quad z, A, B\in \C .\eqno(2.1)$$ Apart of the
classical Lotka-Volterra system, the generalized one splits into
real and complex cases. In appropriate coordinates, the general first
integral of the Lotka-Volterra system in the classical {\it real case} is
$$H(x,y)=x^\lambda y^\mu(1-x-y), \quad \lambda,\mu\in \R,\quad
\lambda\mu(\lambda+\mu+1)\neq 0.\eqno(2.2)$$
The general first integral in the {\it complex case} is
$$H(x,y)=(x^2+y^2)^\lambda
(1-2x)e^{-2\mu {\rm Arctan}(y/x)}, \quad \lambda,\mu\in \R,\quad
\lambda<-\frac12.\eqno(2.3)$$
When $\mu=0$ in the complex case and when
$(\lambda-\mu)(\lambda-1)(\mu-1)=0$ in the real case, the corresponding
system becomes reversible. That is, after a suitable affine change of the
variables, the initial system with a first integral (2.2) or (2.3) will
take in complex coordinates one of the normal forms
$$\dot{z}=-iz+z^2+b\bar{z}^2, \quad b\in \R;\quad \dot{z}=-iz+\bar{z}^2.
\eqno(2.4)$$
This is possible provided that the coefficients in (2.1) satisfy $A^3B\in\R$.
Our main result in this section is the following.

\vspace{2ex} \noindent
 {\bf Theorem 2.} {\it The phase curves of
$(2.1)$ are algebraic curves of genus one if and only if one of the
conditions holds:}
 $$\begin{array}{ll} {\rm (rlv1)} \; A=0 & \mbox{\rm (the Hamiltonian triangle)} \\
{\rm (rlv2)}\; 2AB-\bar{A}^2=0 &   {\rm (lv1)} \; AB+(1\pm2i)\bar{A}^2=0\\
{\rm (rlv3)}\; AB-3\bar{A}^2=0 &   {\rm (lv2)} \; 169AB-(101\pm28i)\bar{A}^2=0\\
{\rm (rlv4)}\; 5AB-3\bar{A}^2=0 &   {\rm (lv3)} \; 289AB-(151\pm42i)\bar{A}^2=0\\
{\rm (rlv5)}\; 5AB-\bar{A}^2=0 &   {\rm (lv4)} \; 1681AB-(783\pm60\sqrt{2}i)\bar{A}^2=0\\
{\rm (rlv6)}\; 3AB+\bar{A}^2=0 &   {\rm (lv5)} \;841AB-(349\pm12i)\bar{A}^2=0 \\
  \end{array}$$
The above statement  is a consequence of the following two propositions
which will be proved together:

\vspace{2ex}
\noindent
{\bf Proposition 1.} {\it The phase curves of the reversible Lotka-Volterra
system $(2.1)$, $A^3B\in\R$, are algebraic curves of genus one if and only if
its first integral is affine equivalent to one of the normal forms}
$$\begin{array}{l} {\rm (rlv1)}\;
H=xy(1-x-y) \qquad \mbox{\rm (the Hamiltonian triangle)} \\
 {\rm (rlv2)}\; H=x^{-3}y(1-x-y) \\
 {\rm (rlv3)}\; H=x^2y(1-x-y) \\
 {\rm (rlv4)}\; H=x^{-4}y(1-x-y) \\
 {\rm (rlv5)}\; H=(x^2+y^2)^{-\frac23}(1-2x)\\
 {\rm (rlv6)}\; H=(x^2+y^2)^{-2}(1-2x).
\end{array}$$

\vspace{2ex}
\noindent
{\bf Proposition 2.} {\it The phase curves of the generic Lotka-Volterra
system $(2.1)$, $A^3B\not\in \R$, are algebraic curves of genus one if
and only if the first integral is affine equivalent  to one of the normal forms}
$$\begin{array}{l}
 {\rm (lv1)}\; H=x^2y^3(1-x-y) \\
{\rm (lv2)}\; H=x^{-6}y^2(1-x-y) \\
 {\rm (lv3)}\; H=x^{-6}y^3(1-x-y) \\
{\rm (lv4)}\; H=x^{-4}y^2(1-x-y) \\
 {\rm (lv5)}\;H=x^{-3}y^\frac32(1-x-y).
\end{array}$$

 \vspace{2ex} \noindent {\bf Proof of Propositions 1 and 2.} Under the conditions in
(2.3), the critical point $\left(\frac{\lambda}{2\lambda+1}, \frac{\mu}{2\lambda+1} \right)$ is a center.  The
origin is a focus for $\mu\neq 0$ and a center elsewhere.  Clearly, the phase curves defined by (2.3) could be
elliptic only if $\mu=0$. Note that if $\mu=0$, the origin is a center for all $\lambda\neq 0$, but this center
is reversible and the system can be transformed into the normal form (1.1). Here we study the Lotka-Volterra
center outside the origin, existing for $\lambda<-\frac12$.

Under the conditions in (2.2) (frankly, one should take
$H=xy(1-x-y)|x|^{\lambda-1}|y|^{\mu-1}$ there,
but modules will be omitted thoroughly), there is a unique critical
point $\left(\frac{\lambda}{\lambda+\mu+1}, \frac{\mu}{\lambda+\mu+1}
\right)$ outside the invariant straight lines $x=0$, $y=0$, $x+y=1$, which
is a center if and only if $\lambda\mu(\lambda+\mu+1)>0$.  In fact, there
are two topologically different configurations having a center, the first
one is obtained for $\lambda>0, \mu> 0$ and the second one corresponds to
the parameters outside the first quadrant.  In this latter case, one can
take without loss of generality $\lambda<0$, $\mu<0$.  Indeed, if e.g.
$\lambda<0<\mu$, then applying an affine change $x=1-X-Y$, $y=Y$, we reduce
the first integral $H(x,y)=t$ in (2.2) to $H_1(X,Y)=t_1$ where
$$H_1(X,Y)=X^{1/\lambda} Y^{\mu/\lambda}(1-X-Y),\quad
t_1=t^{1/\lambda}\eqno(2.5)$$ One can proceed similarly with the other case
$\mu<0<\lambda$.  In the same way, when $\lambda$ and $\mu$ are positive,
we can reduce their values to $\lambda, \mu \in (0,1]$.  Indeed, if e.g.
$\lambda>1$ and $\lambda\geq\mu$, the same change as above transforms (2.2)
into (2.5), with both degrees in (2.5) within $(0,1]$. We proceed similarly if
$\mu>1$ and $\mu\geq\lambda$.

When $\lambda\mu(\lambda+\mu+1)=0$, the first integral in (2.2) should be
replaced by another one containing logarithmic or exponential terms \cite{5,schlomiuk,10}
and therefore its level sets are not algebraic curves.
In this way we have reduced our consideration to the following two cases

 \vspace{2ex}
  (I)  $\;\;0<\lambda\leq 1$, $0<\mu\leq 1$,

\vspace{1ex}
 (II)  $\;\;\lambda<0$, $\mu<0$, $\lambda+\mu+1>0$.

 \vspace{2ex}
 \noindent

In the reversible
Lotka-Volterra cases, one can find all curves with a genus one by simply
using the results from the previous section.  Take the real reversible case
and assume for definiteness that $\lambda=\mu$.  The remaining
possibilities are reduced to this one by the same change which we used to
obtain (2.5) from (2.2).  Let $\lambda=\mu\in (-1/2,1]$.  A further
substitution in (2.2) $x=\frac12(1-X)+Y$, $y=\frac12(1-X)-Y$ transforms the
curve $H(x,y)=t$ to $H_1(X,Y)=t_1$ where $$\textstyle
H_1(X,Y)=X^{1/\lambda}[-\frac{1}{2} Y^2+\frac{1}{8} (X-1)^2], \quad
t_1=\frac{1}{2} t^{1/\lambda}.\eqno(2.6)$$ As
$1/\lambda\in(-\infty,-2)\cup[1,\infty)$ and (2.6) takes the form already
studied in Section 1 above, we conclude immediately that $H_1$, and hence
$H$, is of genus one if and only if $1/\lambda \in\{-4, -3, 1, 2\}$.
Namely, for $\lambda=\mu \in\{1, \frac12, -\frac14, -\frac13\}$ in (2.2).  In
the same way, taking in the complex case (2.3) $\mu=0$ and using the
substitution $x=\frac12(1-X)$, $y=Y$, we reduce (2.3) to (2.6) (with a sign +
in front of $\frac12$).  As $1/\lambda\in(-2,0)$ now, one obtains a curve
of genus one if and only if $1/\lambda \in\{-\frac32, -\frac12\}$.  That
is, for $\lambda \in\{-\frac23, -2\}$ and $\mu=0$ in (2.3).  Thus we have
obtained the complete list of cases with genus one in the reversible
Lotka-Volterra system (the first three entries in the right column contain
the cases with $\mu=1$; the cases with $\lambda=1$ are obtained by a
rotation of the variables $(x,y)\rightarrow (y,x)$)
$$\begin{array}{ll}
H=xy(1-x-y) & \; \\
H=x^\frac12y^\frac12(1-x-y) &  H=x^2y(1-x-y) \\
 H=x^{-\frac14}y^{-\frac14}(1-x-y) & H=x^{-4}y(1-x-y) \\
 H=x^{-\frac13}y^{-\frac13}(1-x-y) & H=x^{-3}y(1-x-y) \\
H=(x^2+y^2)^{-\frac23}(1-2x) &  H=(x^2+y^2)^{-2}(1-2x).
\end{array}$$
For a completeness, let us add to the above list the unique
reversible Lotka-Volterra case having conic orbits, namely
$${\rm (rlv0)}\; H=(x^2+y^2)^{-1}(1-2x).$$
This is the quadratic isochronous center known as ${\cal S}_1$ and it
corresponds to $b=0$ in (2.4). This finishes the proof of Proposition 1.

Now we come to the generic (non-reversible) real cases (I), (II) above, where
$(\lambda-\mu)(\lambda-1)(\mu-1)\neq 0$.  Consider first (I). Let
$$\lambda=\frac{p}{q}, \quad \mu=\frac {r}{s}, \quad p,q,r,s\in \N, \quad
p<q,\quad r<s, \quad gcd(p,q)=1,\quad gcd(r,s)=1.$$
We replace our first integral $H(x,y)=t$ with $H_1(x,y)=t_1$ where:
$$H_1(x,y)=x^{p_1}y^{q_1}(1-x-y)^{r_1},\quad p_1=\lambda r_1, q_1=\mu r_1,
t_1=t^{r_1},\quad r_1=qs/gcd(q,s).\eqno (2.7)$$
The irreducible algebraic curve
$$\Gamma_t= \{(x,y)\in\C^2:x^{p}y^{q}(1-x-y)^{r}=t\},\quad p,q,r \in \N,  \quad p<q<r, \quad gcd(p,q,r)=1
\eqno(2.8)
$$
is smooth if and only if $t\not\in \Delta=\{0, (\frac{r}{p+q+r})^\frac{p^2+q^2+r^2}{r} \}$. Let $\bar{\Gamma}_t$
be the associated compact Riemann surface. Its (geometric) genus is one and the same for all $t\not\in \Delta$.
We have (cf. \cite{3,33}):

\vspace{2ex} \noindent {\bf Proposition 3.}  {\it The compact Riemann surface  $\bar{\Gamma}_t$, $t\not\in
\Delta$, is of genus one if and only if $p=1$, $q=2$, $r=3$.}

\vspace{2ex}
\noindent
{\bf Proof.} For reader's convenience, we repeat the proof from \cite{3}.
It is based on the Poincar\'e-Hopf formula, namely $deg(\omega)=2g-2$,
which we apply to the meromorphic one-form related to the first integral (2.8)
$$\omega= -\frac{dx}{x[q-qx-(q+r)y]}=\frac{dy}{y[p-py-(p+r)x]}.$$
We recall that, in our context,  $deg(\omega)=\sum_k n_k$ where the summation is taken over all zeros or poles
of $\omega$, and $n_k$ are their orders. The above one-form has neither zeroes nor poles in the affine chart
outside the critical locus. Hence it suffices to consider it at infinity. There are 3 points at infinity:
$S_1=[1:0:0]$, $S_2=[0:1:0]$, $S_3=[1:-1:0]$. The local coordinates near $S_1$ can be chosen as follows. Write
$x=1/u$ with $u \rightarrow 0$. After this change of coordinates, equation (2.8) becomes $y^q(1-u+yu)^r=u^{p+r}$
(we take $t=(-1)^r$ for simplicity). Since $u \rightarrow 0$, we have $y \rightarrow 0 $ (if $y \rightarrow
\infty$, we would have $yu \rightarrow -1$ which corresponds  to $S_3$). Consequently, up to bi-analyticity,
there are local coordinates $(u,y)$ verifying $y^{q}=u^{p+r}$. Let $m= gcd(q,p+r)$. There are $m$ different
parameterizations, namely:
$$ u=\xi^{\frac{q}{m}},\quad y=e^\frac{2ik\pi}{m}\xi^\frac{p+r}{m}, \quad
k=0,\ldots, m-1,$$
which correspond to the $m$ different local coordinates near the $m$
smooth points given after blowing-up $S_1$. In these local coordinates,
$$\omega= \frac{-ud(1/u)}{q-q/u-(q+r)y}
=\frac{qd\xi}{m\xi[q-q\xi^{-\frac{q}{m}}
-(q+r)e^\frac{2ik\pi}{m}\xi^{\frac{p+r}{m}}]}.$$ We conclude that
near the $m$ points, $\omega$ has a zero of order $\frac{q}{m}-1$.

We perform the same study near the two other points at infinity.
Near $S_2$, by symmetry, we obtain $n=gcd(p,q+r)$ points where the one-form
$\omega$ has a zero of order $\frac{p}{n}-1$.
Similarly, near $S_3$, there are $l=gcd(r,p+q)$ points where
$\omega$ has a zero of order $\frac{r}{l}-1$.

As a result, we obtain the formula $deg(\omega)=p+q+r-n-m-l$
and therefore the curve (2.8) is elliptic if and only if
$$p+q+r=n+m+l. \eqno (2.9)$$
Now we have to resolve this Diophantine equation.
Obviously, we have $ m \leq q$, $n \leq p$ and $l \leq r$,
hence (2.9) is true if and only if
$$gcd(q,p+r)=q,\quad gcd(p,q+r)=p,\quad gcd(r,p+q)=r.$$
Take natural numbers $\alpha,\beta ,\gamma$ such that
$q+r=p \alpha,\quad  p+r=q \beta, \quad p+q=r \gamma.$
The latter system has a nonzero solution
$$\frac{p}{r}=\frac{\gamma+1}{\alpha+1},\quad
\frac{q}{r}=\frac{\gamma+1}{\beta+1}$$
if and only if
$$\alpha \beta   \gamma =2+\alpha + \beta +\gamma.$$
By (2.8) we have $\alpha>\beta>\gamma$ and the unique solution satisfying this
condition is $\alpha=5$, $\beta=2$, $\gamma=1$ which leads to
$p/r=1/3$, $q/r=2/3$. $\;\square$

\vspace{2ex}
\noindent
To consider the second case (II), we note that a bi-rational
change of the variables $$x=\frac{-X}{1-X-Y},\quad
y=\frac{-Y}{1-X-Y}$$ transforms (II) to a case already considered:
$H_1=X^\Lambda Y^M(1-X-Y)=t_1$ with $\Lambda>0$,
$M>0$, where $$\Lambda=\frac{-\lambda}{\lambda+\mu+1},\quad
M=\frac{-\mu}{\lambda+\mu+1},\quad
t_1=t^{-1/(\lambda+\mu+1)}.\eqno(2.10)$$ Hence one can use the
values of $\Lambda$, $M$ already obtained above to calculate
$\lambda$ and $\mu$ for case (II) through the formula
$$\lambda=\frac{-\Lambda}{\Lambda+M+1},\quad
\mu=\frac{-M}{\Lambda+M+1}.$$ Thus we have completed the list of
all cases with genus one in the generic Lotka-Volterra system (another 15
cases are obtained by a rotation of the variables):
$$\begin{array}{lll}
  H=x^\frac23y^\frac13(1-x-y) & H=x^\frac32y^\frac12(1-x-y)
  & H=x^2y^3(1-x-y) \\
  H=x^{-\frac16}y^{-\frac13}(1-x-y) &
 H=x^{-6}y^2(1-x-y)  & H=x^\frac12y^{-3}(1-x-y) \\
H=x^{-\frac16}y^{-\frac12}(1-x-y) & H=x^{-6}y^3(1-x-y)
& H=x^\frac13y^{-2}(1-x-y) \\
H=x^{-\frac14}y^{-\frac12}(1-x-y) &
H=x^{-4}y^2(1-x-y)  & H=x^\frac12y^{-2}(1-x-y)\\
H=x^{-\frac13}y^{-\frac12}(1-x-y) & H=x^{-3}y^\frac32(1-x-y)
& H=x^\frac23y^{-2}(1-x-y).
\end{array}$$
Proposition 2 is  completely proved.$\;\square$

\vspace{1ex}
\noindent
{\bf Proof of Theorem 2.} The proof is a direct consequence of Propositions 1
and 2.  More precisely, it follows by a straightforward calculation using the
formulas from \cite{5}, case (iv) on page 157 there.  $\;\square$

\vspace{2ex}
 \noindent
 {\large\bf 3. The generating function in the reversible case}

\vspace{2ex} \noindent We are going to study small perturbations
of the reversible system (1.1) in the cases when the first
integral $H$ defines a curve of genus one. Consider first the {\it
quadratic} perturbations of system (1.1) rewritten in real
coordinates: $$\begin{array}{l}\dot{x}=H_y/M+\varepsilon
f(x,y,\varepsilon),\\ \dot{y}=-H_x/M+\varepsilon
g(x,y,\varepsilon)\end{array}\eqno(3.1)$$ where $M=X^{\lambda-1}$
and $f,g$ are quadratic polynomials with coefficients depending
analytically on the small parameter $\varepsilon$. Given a
perturbation $(f,g)$, the limit cycles in (3.1) are determined by
the zeroes of the leading term $I(t)$ in the expansion with
respect to $\varepsilon$ of the displacement map. For this reason
the integral $I(t)$ is called sometimes the {\it generating
function}. As well known \cite{5}, one can always choose a
particular quadratic perturbation so that the corresponding $I(t)$
would have the possible maximum of zeroes within the whole class
of quadratic perturbations. In the generic case $I(t)$ is given by
the complete elliptic integral \cite{5}, Theorems 2 and 3
$$I(t)=\int_{\delta(t)}x^{\lambda-1}(\mu_1+\mu_2
x+\mu_3x^{-1})ydx\eqno(3.2)$$ where $\delta(t)$ is the oval
contained in the level set $H=t$ and $\mu_i\in\R$. In the
exceptional cases (r10) and (r5) with $b=2$ when system (1.1)
belongs to the intersection with the codimension-four stratum
$Q_4$ of the center manifold (see \cite{10} for details), $I(t)$
takes another form and is not an Abelian integral \cite{5}. In the
standard elliptic case when $a=b=-1$ in (1.1) (we will denote it
by (r0)), the integral has a specific form, too. Using (3.2) and
the list of first integrals of genus one (0)-(xviii) above, we
obtain the concrete form of $I(t)$ for all the cases (after
appropriate re-scaling of $x$, $y$, $t$, $H$, $I$) as follows:
$$\begin{array}{ll} {\rm (r0)} &
I(t)=\int_{H=0}(\mu_1+\mu_2t+\mu_3x)ydx,\\ &
H=\frac12y^2-\frac12x^2+\frac13x^3-t,\;\; t\in(-\frac16,0)\\[2mm]
{\rm (r1)} & I(t)=\int_{H=0}x^{-4}(\mu_1+\mu_2
x+\mu_3x^{-1})ydx,\\ {\rm (r2)} &
I(t)=\int_{H=0}x^{-3}(\mu_1+\mu_2x^{-1}+\mu_3x)ydx,\\ & H=\frac12
y^2+\frac{3-b}{6(b+1)}+\frac{b-1}{b+1}x+\frac{1-3b}{2(b+1)}x^2-tx^3\\[2mm]
{\rm (r3)} &
I(t)=\int_{H=0}x^{-4}(\mu_1+\mu_2x^2+\mu_3x^{-2})ydx,\\ {\rm (r4)}
& I(t)=\int_{H=0}x^{-2}(\mu_1+\mu_2x^{-2}+\mu_3x^2)ydx,\\ &
H=\frac12
y^2+\frac{b+3}{24(b+1)}+\frac{b-1}{4(b+1)}x^2+\frac{3b+1}{8(b+1)}x^4-tx^3\\[2mm]
{\rm (r5)} & I(t)=\int_{H=0}x^{-5}(\mu_1+\mu_2x+\mu_3x^{-1})ydx,
\;\; (b\neq 2)\\ {\rm (r5)} &
I(t)=\int_{H=0}x^{-5}[\mu_1+\mu_2x+\mu_3x^{-2}+\mu_4(1-x^{-1})\ln
x]ydx,\;\; (b=2)\\ {\rm (r6)} &
I(t)=\int_{H=0}x^{-4}(\mu_1+\mu_2x^{-1}+\mu_3x)ydx,\\ & H=\frac12
y^2+\frac{2-b}{4(b+1)}+\frac{b-1}{b+1}x+\frac{1-2b}{2(b+1)}x^2-tx^4\\[2mm]
{\rm (r7)} &
I(t)=\int_{H=0}x^{-5}(\mu_1+\mu_2x^3+\mu_3x^{-3})ydx,\\ {\rm (r8)}
& I(t)=\int_{H=0}x^{-2}(\mu_1+\mu_2x^{-3}+\mu_3x^3)ydx,\\ &
H=\frac12
y^2+\frac{1}{12}-\frac13x^3-tx^4,\;\;t\in(-\frac14,0)\\[2mm] {\rm
(r9)} & I(t)=\int_{H=0}(\mu_1+\mu_2x^{-3}+\mu_3x^3)ydx,\\ {\rm
(r10)} &
I(t)=\int_{H=0}x^{-3}[\mu_1+\mu_2x^3+\mu_3x^6+\mu_4(2x^3-1-x^{-3})\ln
x]ydx,\\ & H=\frac12
y^2+\frac16+\frac13x^3-tx^2,\;\;t\in(\frac12,\infty)\\[2mm] {\rm
(r11)} & I(t)=\int_{H=0}x(\mu_1+\mu_2x^{-3}+\mu_3x^3)ydx,\\ {\rm
(r12)} & I(t)=\int_{H=0}x^{-2}(\mu_1+\mu_2x^3+\mu_3x^{-3})ydx,\\ &
H=\frac12
y^2+\frac13+\frac16x^3-tx,\;\;t\in(\frac12,\infty)\\[2mm] {\rm
(r13)} & I(t)=\int_{H=0}(\mu_1+\mu_2x^{-4}+\mu_3x^4)ydx,\\ {\rm
(r14)} & I(t)=\int_{H=0}x^{-4}(\mu_1+\mu_2x^4+\mu_3x^{-4})ydx,\\ &
H=\frac12
y^2+\frac{1}{12}+\frac14x^4-tx^3,\;\;t\in(\frac13,\infty)
\end{array}$$
$$\begin{array}{ll} {\rm (r15)} &
I(t)=\int_{H=0}x^2(\mu_1+\mu_2x^{-4}+\mu_3x^4)ydx,\\ {\rm (r16)} &
I(t)=\int_{H=0}x^{-2}(\mu_1+\mu_2x^4+\mu_3x^{-4})ydx,\\ &
H=\frac12 y^2+\frac14+\frac{1}{12}x^4-tx,\;\;\
t\in(\frac13,\infty) \\[2mm] {\rm (r17)} &
I(t)=\int_{H=0}x^{-5}(\mu_1+\mu_2x^2+\mu_3x^{-2})ydx,\\ {\rm
(r18)} & I(t)=\int_{H=0}x^{-3}(\mu_1+\mu_2x^{-2}+\mu_3x^2)ydx,\\ &
H=\frac12 y^2+\frac16-\frac12x^2-tx^3,\;\; t\in(-\frac13,0).
\end{array}$$
One can formulate the following conjecture concerning the maximal number of
zeroes of the generating function $I(t)$ and the corresponding maximal number
of limit cycles produced by the period annulus (called its {\it cyclicity}).

\vspace{2ex}
\noindent
{\bf Conjecture 1.} (cf. \cite{5}) {\it The period annulus around the center at the
origin in {\rm (r0)-(r22)} has the following cyclicity under small quadratic
perturbations: three for cases {\rm (r1)} with $a^*<a<4$, {\rm (r3)} with
$\frac73<a<4$, {\rm (r4)} with $4<a<5$, {\rm (r5)} with $a=4$, {\rm (r6)}
with $a>4$ and {\rm (r10)}, and two otherwise.}

\vspace{2ex}
\noindent
We note that $a^*\in(\frac53,3)$ is determined from a transcendental equation
\cite{7} and can be calculated numerically, $a^*=2.0655...$

Let us mention that for some of the cases Conjecture 1 has already been
verified. The standard elliptic case (r0) is studied entirely in \cite{11}. The
reversible Hamiltonian case (r2) has been investigated in a series of papers,
see \cite{2} and the references therein. Cases (r5) with $b=2$ and (r10) were
considered in \cite{5} and \cite{6}, respectively. Cases (r5) for $b\neq 2,\frac12$ and
(r6) with $b\in(\frac12,2)$ are studied in \cite{1}. Case (r1) with
$b\in(-1,\frac13)$ and $b\in(\frac13,3)$ are considered in \cite{9} and \cite{7},
respectively. Below we will concentrate our efforts mainly on the cases where
Conjecture 1 remains open and will include in our lists that follow only these
unsolved cases.

\vspace{1ex}
\noindent
In order to reduce the number of Hamiltonians, we replace

$(x,y)\to (1/x, y/x^2)$ in cases (r5) and (r6),

$(x,y)\to (x/(-4t),y/(-4t)^{3/2})$ in cases (r7) and (r8),

$(x,y)\to (2tx,(2t)^{3/2}y)$ in case (r9),

$(x,y,t)\to ((-6s)^{1/3}x,y,(-48s)^{-1/3})$ in cases (r11) and (r12),

$(x,y)\to (3tx,9t^2y)$ in cases (r13) and (r14),

$(x,y)\to (1/(3tx),y/x^2)$ in cases (r15) and (r16),

$(x,y)\to (-x/3t,-y/3t)$ in cases (r17) and (r18).

\noindent
In this way we obtain a reduced list of cases to study as follows:
$$\begin{array}{ll}
{\rm (r1)} & I(t)=\int_{H=0}x^{-4}(\mu_1+\mu_2 x+\mu_3x^{-1})ydx,
\;\;(b\not\in(-1,\frac13)\cup(\frac13,3))\\
{\rm (r11)} & I(t)=\int_{H=0}x(\mu_1+\mu_2t^{-1}x^{-3}+\mu_3tx^3)ydx, \;\;( b=\frac13)\\
{\rm (r12)} & I(t)=\int_{H=0}x^{-2}(\mu_1+\mu_2tx^3+\mu_3t^{-1}x^{-3})ydx,\;\; (b=\frac13)\\
& H=\frac12 y^2+\frac{3-b}{6(b+1)}+\frac{b-1}{b+1}x+\frac{1-3b}{2(b+1)}x^2-tx^3\\[2mm]
{\rm (r3)} & I(t)=\int_{H=0}x^{-4}(\mu_1+\mu_2x^2+\mu_3x^{-2})ydx,\\
{\rm (r4)} & I(t)=\int_{H=0}x^{-2}(\mu_1+\mu_2x^{-2}+\mu_3x^2)ydx,\\
& H=\frac12 y^2+\frac{b+3}{24(b+1)}+\frac{b-1}{4(b+1)}x^2+\frac{3b+1}{8(b+1)}x^4-tx^3\\[2mm]
{\rm (r5)} & I(t)=\int_{H=0}x(\mu_1+\mu_2x^{-1}+\mu_3x)ydx, \;\; (b=\frac12)\\
{\rm (r6)} & I(t)=\int_{H=0}(\mu_1+\mu_2x+\mu_3x^{-1})ydx,\;\;(b\not\in(\frac12,2))\\
& H=\frac12 y^2+\frac{1-2b}{2(b+1)}x^2+\frac{b-1}{b+1}x^3+\frac{2-b}{4(b+1)}x^4-t\\[2mm]
{\rm (r7)} & I(t)=\int_{H=0}x^{-5}(\mu_1+\mu_2t^{-1}x^3+\mu_3tx^{-3})ydx,\\
{\rm (r8)} & I(t)=\int_{H=0}x^{-2}(\mu_1+\mu_2tx^{-3}+\mu_3t^{-1}x^3)ydx,\\
{\rm (r13)} & I(t)=\int_{H=0}(\mu_1+\mu_2tx^{-4}+\mu_3t^{-1}x^4)ydx,\\
{\rm (r14)} & I(t)=\int_{H=0}x^{-4}(\mu_1+\mu_2t^{-1}x^4+\mu_3tx^{-4})ydx,\\
{\rm (r15)} & I(t)=\int_{H=0}x^{-6}(\mu_1+\mu_2t^{-1}x^4+\mu_3tx^{-4})ydx,\\
{\rm (r16)} & I(t)=\int_{H=0}x^{-2}(\mu_1+\mu_2tx^{-4}+\mu_3t^{-1}x^4)ydx,\\
& H=\frac12 y^2-\frac13x^3+\frac14x^4-t,\;\;t\in(-\frac{1}{12},0)\\[2mm]
{\rm (r9)} & I(t)=\int_{H=0}(\mu_1+\mu_2tx^{-3}+\mu_3t^{-1}x^3)ydx,\\
{\rm (r17)} & I(t)=\int_{H=0}x^{-5}(\mu_1+\mu_2t^{-1}x^2+\mu_3tx^{-2})ydx,\\
{\rm (r18)} & I(t)=\int_{H=0}x^{-3}(\mu_1+\mu_2tx^{-2}+\mu_3t^{-1}x^2)ydx,\\
& H=\frac12y^2-\frac12x^2+\frac13x^3-t,\;\; t\in(-\frac16,0).
\end{array}$$
It is useful to know the dimension of the Picard-Fuchs system satisfied by the
basic integrals $I_k(t)=\int_{\delta(t)} x^kydx$ involved in the formulas
above. Let us rewrite the equation $H=0$ into the form
$\frac12y^2={\cal D}(x,t)$. Then for any $k\in\Z$ one obtains
$$\textstyle \int_{\delta(t)} x^k{\cal D}'ydx=\int_{\delta(t)} x^kyd(\frac12y^2)=
\frac13\int_{\delta(t)} x^kdy^3=-\frac{k}{3}\int_{\delta(t)} x^{k-1}y^3dx.$$
Therefore
$$\int_{\delta(t)} x^{k-1}\left(x{\cal D}'+\frac{2k}{3}{\cal D}\right)ydx=0.$$
Taking ${\cal D}=A_0+A_1x+A_2x^2+A_3x^3+A_4x^4$, in terms of the basic
integrals $I_k(t)$ the last identity is equivalent to
$$(2k+12)A_4I_{k+3}+(2k+9)A_3I_{k+2}+(2k+6)A_2I_{k+1}+(2k+3)A_1I_k
+2kA_0I_{k-1}=0. \eqno(3.3)$$
Using (3.3) with $k=0,1,2,\ldots$ and with $k=-1,-2,\ldots$ we obtain that
$$I_k=\mbox{\rm span}(I_0,I_1,I_2),\;\; k\geq 3;\quad
I_k=\mbox{\rm span}(I_{-1}, I_0,I_1,I_2),\;\; k\leq -2 \quad (A_4\neq 0).$$
Similarly, if $A_4=0$, one obtains
$$I_k=\mbox{\rm span}(I_0,I_1),\;\; k\geq 2;\quad
I_k=\mbox{\rm span}(I_{-1}, I_0,I_1),\;\; k\leq -2 \quad (A_4=0, A_3\neq 0),$$
and finally,
$$I_k=\mbox{\rm span}(I_0),\;\; k\geq 1;\quad
I_k=\mbox{\rm span}(I_{-1}, I_0),\;\; k\leq -2 \quad (A_4=A_3=0, A_2\neq 0).$$
Here "span" means in general a polynomial $R[t,t^{-1}]$ span.
As a consequence, we can formulate a result about the dimension of the
Picard-Fuchs system satisfied by the basic integrals in (r1), (r3)-(r9) and
(r11)-(r18).

\vspace{2ex}
\noindent
{\bf Proposition 4.} {\it In cases} (r1), (r3)-(r4) {\it with
$b=-\frac13$}, (r6) {\it with $b=2$}, (r9), (r11)-(r12), (r17)-(r18)
{\it the Picard-Fuchs system is of dimension 3 while in the remaining cases
it is of dimension 4.}

\vspace{2ex}
\noindent
To  derive the Picard-Fuchs equations, we use (3.3) together with the relations
$$I_k'(t)=\int_{\delta(t)}\frac{x^k\partial_t{\cal D}}{y}dx,\quad
\int_{\delta(t)}\frac{x^k{\cal D}'}{y}dx=-kI_{k-1}(t),\quad
\int_{\delta(t)}\frac{x^k{\cal D}}{y}dx=\frac12I_k(t).\eqno(3.4)$$
We can use (3.4) to further simplify some of the integrals above. Thus we
get the final list of reversible cases of genus one yet to study:
$$\begin{array}{ll}
{\rm (r1), (r12)} & I(t)=\int_{H=0}x^{-4}(\mu_1+\mu_2 x^{-1}+\mu_3x)ydx,
\;\;(b\not\in(-1,\frac13)\cup(\frac13,3)) \\
{\rm (r11)} & I(t)=\int_{H=0}x^{-1}(\mu_1+\mu_2x^{-1}+\mu_3x)ydx, \;\;(b=\frac13)\\
& H=\frac12 y^2+\frac{3-b}{6(b+1)}+\frac{b-1}{b+1}x+\frac{1-3b}{2(b+1)}x^2-tx^3,\\[1mm]
&t\in\left(-\frac{2b}{3(b+1)},\frac{2(1-3b)^2}{3(b+1)(3-b)^2}\right), \;\;b<-1;
\quad t\in\left(-\frac{2b}{3(b+1)},0\right), \;\;b\geq \frac13\\[2mm]
{\rm (r3)} & I(t)=\int_{H=0}x^{-4}(\mu_1+\mu_2x^2+\mu_3x^{-2})ydx,\\
{\rm (r4)} & I(t)=\int_{H=0}x^{-2}(\mu_1+\mu_2x^{-2}+\mu_3x^2)ydx,\\
& H=\frac12 y^2+\frac{b+3}{24(b+1)}+\frac{b-1}{4(b+1)}x^2+\frac{3b+1}{8(b+1)}x^4-tx^3,\\[1mm]
& t\in\left(\frac{2b}{3(b+1)},-\frac{2}{3(b+1)}\sqrt{-\frac{3b+1}{b+3}}\right),
\;\;b\in(-3,-\frac13);\\[1mm]
& t\in\left(\frac{2b}{3(b+1)},\infty\right)\;\mbox{\rm otherwise}.\\[2mm]
{\rm (r6)} & I(t)=\int_{H=0}(\mu_1+\mu_2x^{-1}+\mu_3x)ydx,\;\;(b\not\in(\frac12,2))\\
& H=\frac12 y^2+\frac{1-2b}{2(b+1)}x^2+\frac{b-1}{b+1}x^3+\frac{2-b}{4(b+1)}x^4-t,\\[1mm]
&t\in\left(-\frac{b}{4(b+1)},0\right), \;\;b\geq 2;\quad
t\in\left(-\frac{b}{4(b+1)},\frac{(1-2b)^3}{4(b+1)(2-b)^3}\right),
\;\;b\leq \frac12 
\end{array}$$
$$\begin{array}{ll}
{\rm (r7), (r14)} & I(t)=\int_{H=0}x^{-1}(\mu_1+\mu_2x^{-1}+\mu_3x)ydx,\\
{\rm (r8)} & I(t)=\int_{H=0}x^{-1}(\mu_1+\mu_2x^{-1}+\mu_3t^{-1}x^2)ydx,\\
{\rm (r13)} & I(t)=\int_{H=0}(\mu_1+\mu_2x^{-1}+\mu_3t^{-1}x^2)ydx,\\
{\rm (r15)} & I(t)=\int_{H=0}x^{-3}(\mu_1+\mu_2x^{-1}+\mu_3x)ydx,\\
{\rm (r5), (r16)} & I(t)=\int_{H=0}x(\mu_1+\mu_2x^{-1}+\mu_3x)ydx,\\
& H=\frac12 y^2-\frac13x^3+\frac14x^4-t,\;\;t\in(-\frac{1}{12},0)\\[2mm]
{\rm (r9)} & I(t)=\int_{H=0}(\mu_1+\mu_2x^{-1}+\mu_3t^{-1}x)ydx,\\
{\rm (r17)} & I(t)=\int_{H=0}x^{-2}(\mu_1+\mu_2x^{-1}+\mu_3x)ydx,\\
{\rm (r18)} & I(t)=\int_{H=0}(\mu_1+\mu_2x^{-1}+\mu_3x)ydx,\\
& H=\frac12y^2-\frac12x^2+\frac13x^3-t,\;\; t\in(-\frac16,0).
\end{array}$$
We point out that the above formulas are related in each case to the period
annulus around the center at $(1,0)$ obtained for $t\in(t_c, t_s)$ where $t_c$
is the level corresponding to the center and $t_s$ is the level of the contour
at which the annulus terminates.

\vspace{2ex}
 \noindent
 {\large\bf 4. The generating function in the Lotka-Volterra case}

\vspace{2ex}
\noindent
As in the previous section, we will obtain the complete list of generating
functions $I(t)$ for the generalized Lotka-Volterra systems (2.1) in the cases
when all their orbits are elliptic curves.  Apart of the reversible case, the
equivalence classes for the Lotka-Volterra case listed in Propositions 1 and 2
are a result of {\it affine} transformations.  Therefore we can choose by one
representative from each such a class and then use the generating function
corresponding to it.  We will always choose the cases whose parameters satisfy
conditions (I) and (II) above.  Take real coordinates $z=x+iy$ and consider a
small quadratic perturbation as in (3.1).  In the reversible case (2.4), the
first integral and integrating factor are respectively \cite{5}
$$H(X,y)=X^\lambda\left(\frac{y^2}{2}-\frac{\lambda(\lambda-1)^2}{32(\lambda+2)}
\left(X-\frac{\lambda+2}{\lambda}\right)^2\right),\quad M=X^{\lambda-1}
\eqno(4.1)$$ where $$\lambda=\frac{b+1}{b-1},\quad b\neq 0, \pm1,\frac13,\quad
X=1-\frac{4x}{\lambda-1}.$$ The generating function for small quadratic
perturbations reads \cite{5}, Theorem 3:
$$I(t)=\int_{H=t}x^{\lambda-1}[\mu_1y+\mu_2yx^{-1}+\mu_3(x-1)y^{-1}]dx.
\eqno(4.2)$$ In the Hamiltonian triangle case (rlv1), the integral $I(t)$
takes another form and this case has already been studied in \cite{14}.  For
this reason (rlv1) is not included in our list below.  By (4.1), (4.2), (2.6)
and the formulas of first integrals of genus one (rlv2)--(rlv6), we get the
list of explicit expressions of $I(t)$ for all other reversible Lotka-Volterra
cases as follows: $$\begin{array}{ll} {\rm (rlv2)} &
I(t)=\int_{H=t}x^{-4}(\mu_1y+\mu_2yx^{-1}+\mu_3(x-1)y^{-1})dx,\\ &
H=x^{-3}(\frac12y^2-(x-\frac13)^2),\;\; t\in(-\frac49,0)\\[2mm] {\rm (rlv3)} &
I(t)=\int_{H=t}x(\mu_1y+\mu_2yx^{-1}+\mu_3(x-1)y^{-1})dx,\\ &
H=x^2(\frac12y^2-(x-2)^2),\;\; t\in(-1,0)\\[2mm] {\rm (rlv4)} &
I(t)=\int_{H=t}x^{-5}(\mu_1y+\mu_2yx^{-1}+\mu_3(x-1)y^{-1})dx,\\ &
H=x^{-4}(\frac12y^2-(x-\frac12)^2),\;\; t\in(-\frac14,0)\\[2mm] {\rm (rlv5)} &
I(t)=\int_{H=t}x^{-\frac52}(\mu_1y+\mu_2yx^{-1}+\mu_3(x-1)y^{-1})dx,\\ &
H=x^{-\frac32}(\frac12y^2+(x+\frac13)^2),\;\; t\in(\frac{16}{9},\infty)\\[2mm]
{\rm (rlv6)} &
I(t)=\int_{H=t}x^{-\frac32}(\mu_1y+\mu_2yx^{-1}+\mu_3(x-1)y^{-1})dx,\\ &
H=x^{-\frac12}(\frac12y^2+(x+3)^2),\;\; t\in(16,\infty).  \end{array}$$
Clearly the elliptic curves in cases (rlv5) and (rlv6) are obtained by
introducing a new variable $x\to x^2$ so that they are given by the level sets
of $H(x^2,y)$.  Let us also recall that the functions $I(t)$ in (rlv2)-(rlv6)
are the coefficients at $\varepsilon^2$ in the expansion of the displacement
map obtained for the special perturbations in (3.1) which keep the coefficient
at $\varepsilon$ zero.

\vspace{1ex}
\noindent
What concerns the small quadratic perturbations of the generic (nonreversible)
Lotka-Volterra system, in (3.1) we have $M=x^{\lambda-1}y^{\mu-1}$
(modules needed outside the first quadrant but we will omit them) where
$H$ is determined by (2.2). In the generic case $I(t)$ is given by the complete
elliptic integral \cite{5}, Theorem 2
$$I(t)=\int\!\!\int_{Int\;\delta(t)}x^{\lambda-1}y^{\mu-1}(\mu_1+\mu_2 x^{-1}
+\mu_3y^{-1})dxdy.\eqno(4.3)$$
Using (4.3) and the list of first integrals of genus one (lv1)--(lv5),
we easily obtain the concrete form of $I(t)$ for all the cases as follows:
$$\begin{array}{ll}
{\rm (lv1)} & I(t)=\int_{H=t}x^{-\frac13}y^{\frac13}(\mu_1+\mu_2x^{-1}+\mu_3y^{-1})dx,\\
 & H=x^\frac23y^\frac13(1-x-y), \;\;t\in(0,432^{-1/3}),\\[2mm]
{\rm (lv2)} & I(t)=\int_{H=t}x^{-\frac76}y^{-\frac13}(\mu_1+\mu_2x^{-1}+\mu_3y^{-1})dx,\\
 & H=x^{-\frac16}y^{-\frac13}(1-x-y),  \;\;t\in(432^{1/6},\infty),\\[2mm]
{\rm (lv3)} & I(t)=\int_{H=t}x^{-\frac76}y^{-\frac12}(\mu_1+\mu_2x^{-1}+\mu_3y^{-1})dx,\\
 &  H=x^{-\frac16}y^{-\frac12}(1-x-y), \;\;t\in(432^{1/6},\infty), \\[2mm]
{\rm (lv4)} & I(t)=\int_{H=t}x^{-\frac54}y^{-\frac12}(\mu_1+\mu_2x^{-1}+\mu_3y^{-1})dx,\\
 & H=x^{-\frac14}y^{-\frac12}(1-x-y),  \;\;t\in(2\sqrt2,\infty),\\[2mm]
{\rm (lv5)} & I(t)=\int_{H=t}x^{-\frac43}y^{-\frac12}(\mu_1+\mu_2x^{-1}+\mu_3y^{-1})dx,\\
 & H=x^{-\frac13}y^{-\frac12}(1-x-y), \;\;t\in(432^{1/6},\infty).  \end{array}$$
As in Section 3, we formulate a general conjecture about the cyclicity
of the period annulus in the Lotka-Volterra systems having all their
orbits elliptic or conic curves.

\vspace{2ex}
\noindent
{\bf Conjecture 2.}  {\it The cyclicity under small quadratic perturbations
of the period annulus surrounding the center at the origin in the
Lotka-Volterra system $(2.1)$ is as follows: three in case} (rlv1) {\it $($the
Hamiltonian triangle$)$ and two in all other cases} (rlv2)-(rlv6),
(lv1)-(lv5).

\vspace{2ex}
\noindent
Except for the isochronous center ${\cal S}_1$ and the Hamiltonian triangle,
it is very likely that Conjecture 2 is still open.  We recall that quadratic
perturbations of the general Lotka-Volterra system have been considered in
\cite{10}.  However, in the recent book \cite{13}, page 379, V.I.  Arnold
declared that the problem with the Lotka-Volterra system still remains open.

\vspace{2ex}
\noindent
{\bf Remark about the codimension four case.} In this remark we discuss
in brief the generic (non-reversible) codimension 4 center.  In complex
coordinate $z=x+iy$, the related system becomes
$$\dot{z}=-iz+4z^2+2|z|^2+\alpha\bar{z}^2,\quad \alpha\in\C\setminus \R,
\quad |\alpha|=2.$$
In $(\bar{x},\bar{y})=(X^2,Y)$ coordinates, where
$$Y=(2+b)x+cy,\quad X=1+8x+\frac{4}{2+b}Y^2, \quad \alpha=b+ic,$$
the system has a first integral of the form \cite{5} (below the bars are
omitted)
$$H(x,y)=\frac{x^{-3}}{8(2-b)}\left(\frac{4y^3}{3(2+b)}+\frac{4y^2}{2+b}
+(1-x^2)y-x^2+\frac13\right).$$
It is seen that the level sets of this first integral are (generically)
elliptic curves. Therefore the generating function $I(t)$ whose zeroes
correspond to limit cycles in the perturbed system are given by the
following complete elliptic integral (cf. \cite{5}, Theorem 2 (iii))
$$I(t)=\int\!\!\int_{H(x,y)<t}x^{-6}[\mu_1+\mu_2y+\mu_3y^3+
\mu_4(\kappa^2y^4-x^4)]dxdy, \;\;\kappa=\frac{4}{2+b}.$$ The conjecture that the integral $I(t)$ has at most
three zeroes in the interval $(t_c,t_s)=(-\frac{1}{12(2-b)},0)$ corresponding to the period annulus around the
center at $(x_c,y_c)=(1,0)$ is still open.

\vspace{2ex}
{\large \bf 5.  Zeros of Abelian integrals for some of the reversible cases}

\vspace{2ex} \noindent
In this section we study the generating functions $I(t)$ related to the
reversible quadratic systems (r18) and (r11). These systems contain no
parameter, have a unique period annulus around the center at $(1,0)$ and
the Picard-Fuchs system satisfied by the components of $J(t)=I'(t)$ is of
dimension three and two (respectively). Our main result  is the following.

\vspace{2ex} \noindent {\bf Theorem 3.} {\it The exact upper bound of
the number of the limit cycles produced by
the period annulus under quadratic perturbations of the reversible
system $\rm(r18)$ or $\rm(r11)$ is two.}

\vspace{2ex} \noindent
To prove Theorem 3 we study first one-parameter analytic quadratic
perturbations of (r18) and (r11). According
to the formulas derived in Section 3, the number of the limit cycles
of such perturbations is bounded by the
number of the zeros of the following generating functions

$$I(t)=\int_{\delta(t)}(\mu_1+\mu_2x^{-1}+\mu_3x)ydx,\;\;
t\in \left(-\frac16,0\right), \;\; H_t=\frac{1}{2}y^{2}-\frac{1}{2}x^{2}+\frac{1}{3}x^{3}-t, \eqno(5.1)$$
for system (r18) and
$$I(t)=\int_{\delta(t)}x^{-1}(\mu_1+\mu_2x^{-1}+\mu_3x)ydx,\;\;
t\in \left(-\frac16,0\right),\;\; H_t=\frac{1}{2}y^{2}+\frac{1}{3}-\frac{1}{2}x-tx^{3},  \eqno(5.2)$$
for system (r11). In the above integrals $\delta(t)$ denotes the unique oval
of the real algebraic curve $\{(x,y)\in\R^2:H_t(x,y)=0\}$ for a given
$t\in(-\frac16,0)$ and $\mu_i$ are arbitrary real constants.
The ovals $\delta(t)$ in (5.1) form a period annulus which is bounded by the
homoclinic loop $\{H_0=0\}$ going through the saddle at the origin,
while in the case (5.2) they form a period annulus bounded by the
parabola $\{H_0=0\}$. In both cases, the ovals $\delta(t)$ exist for
$t\in(-\frac16,0)$, $I(t)$ is analytic in
a neighborhood of $t=-1/6$, and $I(-1/6)=0$. Consider the derivative
$J(t)=I'(t)$, see formulae (5.12), (5.13)
bellow. The key ingredient of the proof of Theorem 3 is the following result.

\vspace{2ex} \noindent
{\bf Theorem 4.} {\it The three-dimensional vector space of Abelian integrals
$J(t)=I'(t)$, $t \in[-\frac16,0)$, defined by $(5.1)$ or $(5.2)$ is Chebyshev.
This means that each integral $J(t)$ has at most two zeros $($counted with
multiplicity$)$ in the interval $[-\frac16,0)$.}

\vspace{2ex} \noindent
Before proving Theorem 4 we need some preparation.

\vspace{2ex} \noindent
{\bf 5.1.  Monodromy of the level curves.} Let us denote
$$\Gamma_t=\{(x,y) \in \C^{2},\; H_t(x,y)=0\}, \quad
 \overline{\Gamma}_t=\{[x:y:z] \in \C\P^2,\;
H_t\left(\frac{x}{z},\frac{y}{z}\right)=0\}.$$
If $t\neq -\frac16, 0$, then both $\Gamma_t$ and $\overline{\Gamma}_t$ are
smooth curves and $\overline{\Gamma}_t$ is a compact Riemann surface of genus
one which is also the compactification of the affine elliptic curve $\Gamma_t$.
We have $\overline{\Gamma}_t=\Gamma_t \cup \infty$ where $\infty=[0:1:0]$ and
therefore
$$rank\; H_1(\Gamma_t, \Z)=2, \quad rank\; H_1(\overline{\Gamma}_t, \Z)=2.$$
For (5.2), all compact complex level curves $\overline{\Gamma}_t$ have two
complex conjugate common points, namely: $P^{+}=(0, i \sqrt{\frac{2}{3}})$ and
$P^{-}=(0, -i \sqrt{\frac{2}{3}})$ in affine coordinates, as well as one common
point at infinity. For (5.1) there is one common point, the point at infinity.
For both cases, when blowing-up these points on the complex projective plane
$\P^2$, one obtains a compact smooth rational surface $S$ and an analytic map
$ S \stackrel{\pi}{\rightarrow} \P $ where the projection $\pi$ is induced by
one of the rational maps
$$ \C^2 \rightarrow \C: (x,y) \mapsto t=
\frac{1}{2}y^{2}-\frac{1}{2}x^{2}+\frac{1}{3}x^{3}, \eqno(5.3) $$
$$ \C^2 \dashrightarrow \C: (x,y) \mapsto t=
\frac{\frac{1}{2}y^{2}+\frac{1}{3}-\frac{1}{2}x}{x^{3}}.  \eqno(5.4) $$
In both cases $S$ is an elliptic surface in the sense of Kodaira
\cite{singularbis}
with three singular fibers $\pi^{-1}(0), \pi^{-1}(-1/6), \pi^{-1}(\infty)$.
In particular, this allows one to compute the global monodromy of the related
homology bundle as described below.

\vspace{2ex} \noindent {\bf 5.1.1. The monodromy of the fibration (5.3).} Denote by $0^\pm$ the two points on
$\Gamma_t$ with coordinates $(x=0,y=\pm\sqrt{2t})$ and let $\widehat{\Gamma}_t=\overline{\Gamma}_t \setminus
\{0^+,0^-\}$. We shall determine the monodromy of the homology bundle associated to $(5.3)$ with fibers
$H_1(\widehat{\Gamma}_t,\Z) = \Z^3$.  Let $\delta(t)$, $\gamma(t)$ be a continuous family of cycles, vanishing
at $t=-1/6$ and $t=0$ respectively, and let $\alpha(t)$ be a cycle represented by a small loop around $0^+$ on
$\widehat{\Gamma}_t$.  We need to describe the monodromy of these cycles on the plane $\C\setminus\{0,-1/6\}$.

The precise definition of the families  $\delta(t)$, $\gamma(t)$ is as follows. For $t\in(-1/6,0)$ the
polynomial $\frac{1}{3}x^{3}-\frac{1}{2}x^{2}-t$ has three real roots $x_1(t) < 0 < x_2(t) < x_3(t)$. Let $l$ be
a simple loop on $\C\setminus\{x_1(t), 0, x_2(t), x_3(t)\}$ which makes one turn about $x_1(t), 0, x_2(t)$ in a
positive direction and does not contain in its interior $x_3(t)$.

We define $\gamma(t)$ to be the cycle represented by $l$ on
$\widehat{\Gamma}_t$ (that is to say by one of the two pre-images of $l$ under
the projection $(x,y)\mapsto x$).  The cycle $\gamma(t)$ is defined up to a
sign.  The cycle $\pm \delta(t)$ is defined in a similar way, but the simple
loop $l$ is supposed to make one turn about $x_2(t), x_3(t)$ in a positive
direction, and does not contain $x_1(t), 0$ in its interior.

We note that $x_1(0)=x_2(0)=0$, $x_3(0)>0$ and hence $\gamma(t)$ is a cycle
vanishing at $t=0$.  Similarly, $\delta(t)$ is a vanishing cycle at $t=-1/6$,
as $x_1(-1/6)<0$, $x_2(-1/6)=x_3(-1/6)>0$.

 Let
$$\Pi:\; \pi_1(\C\P^1 \setminus \Delta,t_0) \rightarrow
Aut\,(H_1(\widehat{\Gamma}_t, \Z)),\;\; \textstyle
\Delta=\{-\frac16, 0, \infty\},\; t_0\in (-\frac16,0)$$ be the monodromy
representation related to the elliptic
fibration associated to $(5.3)$. The image $$\Pi( \pi_1(\C\P^1
\setminus \Delta,t_0))$$  is a group
generated by $l_0^{*}, l_1^{*}, l_{\infty}^{*}= l_1^*\circ l_0^*$, the
monodromy operators corresponding respectively to the oriented loops $l_0, l_1,
l_{\infty}= l_1\circ l_0$.  Here $l_0$ is a simple loop which makes one turn
about $0$ in a positive direction and does not contain $-1/6$ in its interior.
Similarly $l_1$ is a simple loop which makes one turn about $-1/6$ in a
positive direction and does not contain $0$ in its interior.

\vspace{2ex} \noindent
{\bf Lemma 1.} {\it In the global basis $(\alpha(t),\delta(t), \gamma(t))$
$($with appropriate orientation of the cycles$)$, the monodromy operators
associated to the fibration defined by $(5.3)$ are: }
 $$l_0^{*}: M_0=\left ( \begin{array}{ccc} -1& 1 & 0 \\ 0 &1 & 0\\ 0 &-1&1
 \end{array}  \right ),\quad
 l_1^{*}: M_1=\left ( \begin{array}{ccc} 1& 0 & 0 \\ 0 & 1 & 1\\ 0 & 0 & 1
 \end{array}  \right ), $$
 $$ l_{\infty}^{*}: M_\infty= M_1 M_0 =
\left(\begin{array}{ccc} -1& 1 & 1 \\ 0 &0 & 1\\ 0&-1& 1\end{array}\right).$$

\vspace{2ex} \noindent
{\bf Proof.} The identity
$$2\pi i \int _{\alpha(t)}\frac{dx}{x y}=
\frac{1}{\sqrt t}$$ shows that $l_0^{*} (\alpha(t)) = -\alpha(t)$.  The
remaining claims follow from the usual Picard-Lefschetz formula (it is enough
to describe the monodromy of the roots $x_i(t)$) and therefore the details
are omitted. $\Box$

\vspace{2ex} \noindent
{\bf 5.1.2. The monodromy of the fibration (5.4).}
Here we determine the monodromy representation
$$\Pi: \pi_1(\C\P^1 \setminus \Delta,t_0) \rightarrow
Aut\,(H_1(\overline{\Gamma}_t, \Z)), \textstyle
\;\;\Delta=\{-\frac16, 0, \infty\}, t_0 \in (-\frac16,0)$$
associated to the elliptic fibration defined by $(5.4)$.
The automorphisms $l_0^{*}, l_1^{*}, l_{\infty}^{*}$ associated to the
oriented loops $l_0, l_1, l_{\infty}$ are defined as in the preceding
subsection.

For $(5.4)$ the global monodromy does not follow from the Picard-Lefschetz
formula.  The local monodromy (around the singular fibers), however, depends
only on the topological type of these fibers and is computed in
\cite{singularbis}. The topological type of the fibers $\pi^{-1}(0)$,
$\pi^{-1}(-1/6)$, $\pi^{-1}(\infty)$ on its turn is computed in \cite{3} and
it is  respectively (III), $\rm (I_1)$, $\rm(IV^*)$ (using the Kodaira notations,
see e.g. \cite[Table 6.1]{miranda}). We conclude that, up to a conjugation,
the monodromy operators are given by (e.g.  \cite[Table 1]{H})
$$\begin{array}{ll}
l_\infty^*: M_\infty=\left(\begin{array}{cc}-1 & -1 \\ 1 & 0\end{array}\right)
& \rm(IV^*), \\[4mm] l_1^*: M_1=\left(\begin{array}{cc}1 & 1 \\ 0 &
1\end{array}\right) & \rm(I_1),\\[4mm] l_0^*: M_0=\left(\begin{array}{cc} 0 & 1
\\ -1 & 0 \end{array}\right) & \rm(III).\end{array}$$
By abuse of notation, we
denote by $\delta(t) \in H_1(\overline{\Gamma}_t, \Z)$, $t\in
(-1/6,0)$, the continuous family of cycles, vanishing when $t \rightarrow
-\frac16$ along a path connecting $t$ to $-\frac{1}{6}$ in $(-\frac16,0)$. Let
$\gamma(t)=l_0^{*}(\delta(t))$. The continuous families of cycles $\delta(t),
\gamma(t)$ are defined in $\C\setminus \{-\frac16,0\}$.  According to
\cite{movasati}, $(\delta(t), \gamma(t))$ is a basis of
$H_1(\overline{\Gamma}_t, \Z)$.

\vspace{2ex} \noindent
{\bf Lemma 2.} {\it In the global basis $(\delta(t), \gamma(t))$ $($with
appropriate orientations of the cycles$)$, the monodromy operators associated
to the fibration defined by $(5.4)$ are:}

 $$l_0^{*}:M_0=\left(\begin{array}{cc} 0 & -1 \\1 & 0 \end{array}\right),\;
 l_1^{*}:M_1=\left (\begin{array}{cc} 1 & -1 \\0 & 1 \end{array}\right),\;
 l_{\infty}^{*}:M_\infty= l_1^{*} \circ l_0^{*} =
 \left ( \begin{array}{cc} -1 & -1 \\ 1 & 0 \end{array}\right).$$

\vspace{2ex} \noindent
{\bf Proof.} For the global sections $\delta, \gamma$ of the homology bundle we
have $l_0^*\delta =\gamma$, $l_1^*\delta = \delta$ and hence
$$
M_0=\left ( \begin{array}{cc} 0 & * \\1 & * \end{array}  \right ),\quad
M_1=\left ( \begin{array}{cc} 1 & * \\0 & * \end{array}  \right ).
$$
The relations $Tr(M_1)=2$, $Tr(M_0)= 0$, $det(M_0)=1$,
$Tr(M_1M_0)=Tr(M_\infty)=-1$ imply the result.$\;{\square}$\\

\vspace{2ex} \noindent {\bf 5.2. Wronskians.} Let
$$\omega_k=\frac{x^{k}  dx}{y} ,\;\; k \in \Z$$
be polynomial one-forms on $\C^2$.  They induce meromorphic differential
one-forms on the compact Riemann surface $\overline{\Gamma}_t$ which are
denoted in the same way.

Let us introduce the  Wronskian
$$W_{\delta(t), \gamma(t)}(\omega_k,\omega_0)=
\int_{\delta(t)}\omega_k \int_{\gamma(t)}\omega_0 -\int_{\delta(t)}\omega_0
\int_{\gamma(t)}\omega_k$$ where, by
abuse of notation, $\delta(t), \gamma(t)$ are continuous families of
\emph{closed loops}, with the properties :
\begin{itemize}
    \item they intersect transversally at a single point
    \item $\overline{\Gamma}_t\setminus \{\delta(t), \gamma(t)\}$ is
    homeomorphic to a rectangle whose opposite sides are identified to
    $\delta(t)$ and $ \gamma(t)$ \item The poles of $\omega_k$ are contained
    in the interior of this rectangle
\end{itemize}
As $\omega_0$ is holomorphic in both cases, the above Wronskian can be
computed by making use of the reciprocity law for abelian differentials of the
second and third kind \cite[p.  647]{4} as shown next.  Below, we denote by $C$
a certain nonzero constant.  In the case (5.2) the Wronskian is a rational
function and
$$
W_{\delta, \gamma}(\omega_k,\omega_0)=
\frac{Res_\infty \omega_k \int_\infty^P \omega_0}{2\pi i } =
\left \{\begin{array}{ll}
\frac{C}{t},  & k=1, \\
0, & k=2, \\
\frac{C}{t^2},  & k=3. \\
\end{array}\right.
\eqno(5.5)
$$
In a similar way, in the case (5.1) we have
$$
W_{\delta, \gamma}(\omega_{-1},\omega_0)= \frac{Res_{0^-} \omega_{-1}}{2\pi i } \int_{P_0}^{0^-} \omega_0 +
\frac{Res_{0^+} \omega_{-1} }{2\pi i }\int_{P_0}^{0^+} \omega_0 = \frac{1}{\sqrt{2t}} \int_{0^-}^{0^+} \omega_0
\eqno(5.6)
$$
and
$$
W_{\delta, \gamma}(\omega_{1},\omega_0) = \frac{Res_\infty \omega_1 \int_\infty^P \omega_0}{2\pi i } = C
\neq 0 . \eqno(5.7)
$$

\vspace{2ex} \noindent {\bf 5.3. Asymptotics of the Abelian integrals.} Here we study, for suitable $k$, the
asymptotical behavior of $J_k(t)=\int_{\delta(t)}\omega_k,$ where $\delta(t)$ is the continuous family of cycles
vanishing at the center (associated to $t=-1/6$), near $t= \infty$ and $t=0$.  Below, we shall write $J(t)
\lesssim (t-t_0)^{s}\log(t-t_0)^{r}, r, s \in \mathbb{R}$ provided that for every sector $S$ centered at the
critical value $t_0 \in \Delta$ there exists a nonzero constant $C_S$ such that $ |J(t)| \leq C_S
|t-t_0|^{s}|\log(t-t_0)^{r}|$.

\vspace{2ex} \noindent
{\bf 5.3.1. The case (5.1).} Here we have
$$J_k(t)=\int_{\delta(t)}\frac{x^{k}dx}{y}=\int_{\delta(t)}\frac{x^{k}dx}{
\sqrt{-\frac{2}{3}x^{3}+x^{2}+2t}}, \,\,\,\, k=0,\pm1 .$$

\noindent
\begin{description}
\item[(a)]
 Near $t= \infty$, the change of variables $x=t^{\frac{1}{3}}u$ leads to
$$J_k(t)=t^{\frac{k}{3}-\frac{1}{6}}\int_{\tilde{\delta}(t)}\frac{u^{k}du}{
 \sqrt{-\frac{2}{3}u^{3}+t^{-\frac{1}{3}}u^{2}+ 2}}=t^{\frac{k}{3}
 -\frac{1}{6}}\tilde{J}_k(t).$$
When $|t| \rightarrow \infty$ (with bounded argument) the integral $\tilde{J}_k(t)$ tends to a finite constant.
Consequently
$$ J_k(t) \lesssim t^{\frac{k}{3}-\frac{1}{6}},\quad k \in \Z. \eqno(5.8)$$

\item[(b)] Near $t=0$,  the Abelian integral $J_k(t)$  can be expanded  as follows
$$
J_k(t)= -\frac{\ln (t)}{2\pi i} \int_{\gamma(t)} \omega_k - \frac12 \int_{\alpha(t)}\omega_k + Q(t)
$$
where $Q$ is a meromorphic function in a neighborhood of $t=0$ (this follows from Lemma 1). Therefore
$$
J_k(t)= \frac{C}{\sqrt{t}} + P(t) \log(t) + Q(t)
$$
where $P,Q$ are meromorphic in a neighborhood of $t=0$. We claim that
$$J_k(t) \lesssim \left \lbrace \begin{array}{ll}
\log(t) &\ \textrm{if} \  k = 0,\\
1 &\ \textrm{if} \  k = 1,\\
 \frac{1}{\sqrt{t}} &\ \textrm{if} \ k=-1.
\end{array}\right. \eqno(5.9)$$
Indeed, for $k\geq 1$
$$
\lim_{t\rightarrow 0^-} J_k(t)= 2 \int_0^{x_3(0)} \frac{x^{k-1}dx}{
\sqrt{-\frac{2}{3}x+1}} $$
where $x_3(0)=3/2$.  For $k=0$ the integral $J_0(t)$ is of the first kind and
its behavior is well known.  For $k=-1$ the leading term of the expansion of
$J_{-1}(t)$ is given by
$$ \int_{\sqrt{t}}^\infty \frac{dx}{x\sqrt{x^2-t}}=
\frac{1}{\sqrt{t}} \int_{1}^\infty \frac{dx}{x\sqrt{x^2-1}}.  $$

\item[(c)] Near $t=-1/6$ the functions $J_k(t)$ are holomorphic and $J_k(t) \lesssim 1$.
\end{description}
\vspace{2ex} \noindent
{\bf 5.3.2. The case (5.2).}  We have
$$J_k(t)=\int_{\delta(t)}\frac{x^kdx}{y}=\int_{\delta(t)}\frac{x^kdx}{
\sqrt{2tx^{3}+x-\frac{2}{3}}},\quad k=1,2,3.$$
\begin{description}
\item[(a)]
Near $t= \infty$, the change of variables $x=t^{-\frac{1}{3}}u$ leads to
$$J_k(t)=t^{-\frac{k+1}{3}}\int_{\tilde{\delta}(t)}\frac{u^{k}du}{
\sqrt{2u^{3}+t^{-\frac{1}{3}}u-\frac{2}{3}}} =t^{-\frac{k+1}{3}}\tilde{J}_k(t).$$ When $|t| \rightarrow \infty$
(with a bounded argument) the integral $\tilde{J}_k(t)$ tends to a finite constant (possibly zero) and hence
$$
J_k(t) \lesssim t^{-\frac{k+1}{3}}, \;\; k \in \Z. \eqno(5.10)
$$
\item[(b)]
 Near $t=0$ any Abelian integral of the first or second kind is a power series in $t^{1/4}$ (because the
eigenvalues of $l_0^*$ are fourth roots of the unity). As $t$ tends to $0$ two roots of $tx^3+x/2-1/3$ tend to
$\infty$ which suggests us to make the change of variables $x=t^{-\frac{1}{2}}u$. Then we have
$$J_k(t)=t^{-(\frac{k}{2}+\frac{1}{4})}\int_{\tilde{\delta}(t)}\frac{u^{k}du}{
\sqrt{2u^{3}+u-\frac{2}{3}t^{\frac{1}{2}}}} =t^{-(\frac{k}{2}+\frac{1}{4})}\tilde{J}_k(t)$$  and hence
$$
J_k(t) \lesssim  t^{-(\frac{k}{2}+\frac{1}{4})}, \;\;\;  k \in \Z.
\eqno(5.11) $$
\item[(c)] Near $t=-1/6$ the integrals $J_k(t)$ are holomorphic and $J_k(t) \lesssim 1$.
\end{description}

\vspace{2ex} \noindent {\bf 5.4. Proof of Theorem  4.} We have to prove that
the derivative
$$
J(t)=I'(t)= \mu_1 \int_{\delta(t)}  \frac{ x^2 dx}{y}+ \mu_2 \int_{\delta(t)}
\frac{ x dx}{y} + \mu_3 \int_{\delta(t)}  \frac{ x^3 dx}{y}  \eqno{(5.12)}
$$
(in the case (r11)), or
$$
J(t)=I'(t)= \mu_1 \int_{\delta(t)}  \frac{  dx}{y}+ \mu_2 \int_{\delta(t)}
\frac{ dx}{x y} + \mu_3 \int_{\delta(t)} \frac{ x dx}{y} \eqno{(5.13)}
$$
(in the case (r18)) has at most two zeros. For this, we use the  method of Petrov, see e.g. \cite{8}, based on
the argument principle.

Introduce the function
$$
F(t)= \frac{J(t)}{J_0(t)},\quad J_0(t)= \int_{\delta(t)}  \frac{  dx}{y}.
$$
It is real analytic on $(-1/6,0)$ and has an analytic continuation in the complex domain
$\mathcal{D}=\C\setminus[0,\infty)$, because the period $J_0(t)$ of the elliptic curve $\bar{\Gamma}_t$ does not
vanish there, including at the point $t=-1/6$. To bound the number of the zeros of $F$ in $\mathcal{D}$ it is
enough to find this number in the smaller domain $\mathcal{D}_{Rr}= \mathcal{D}\cap \{t: r< |t| < R\}$ for $r$
sufficiently small and $R$ sufficiently big. We are going to evaluate the increment of the argument of $F$ along
the boundary of $\mathcal{D}_{Rr}$, oriented in a positive direction.

\vspace{2ex} \noindent {\bf 5.4.1. Zeros of $F$ in the case (r18) and (5.1).} Along the boundary of the small
circle $\{|t|=r\}$, according to $(5.9)$, the increase of the argument of $F$ is at worst close to $\pi$, and
along the boundary of the big circle $\{|t|=R\}$ this increase is at worst close to $2\pi/3$, see $(5.8)$.
Denote by $F^\pm$ the restriction on $(0,\infty)$ of the analytic function obtained as an analytic continuation
of $F$ along an arc contained in the upper (lower) complex half-plane ${\rm Im}\, \pm t > 0$ respectively.
 The function $F$ is real analytic on $(-\infty,0)$ which implies
that along the interval $(0,\infty)$ the imaginary part ${\rm Im}\,(F^+(t))$ of $F^+$ satisfies
$$ 2i\,{\rm Im}\,(F^+(t))=F^+(t)-\overline{F^+(t)}= F^+(t)-F^-(t).$$
Assume for a moment that $F^+(t)-F^-(t)$ has at most one simple zero on the interval $(0,\infty)$. Then
\emph{summing up the above information we conclude that $F(t)$ has at most two zeros in the complex domain
$\mathcal{D}_{Rr}$, and hence in $\mathcal{D}$}. This result is obviously exact.

It remains to prove that $F^+(t)-F^-(t)$ has at most one zero on $(0,\infty)$. Clearly $F^-$ is an analytic
continuation of $F^+$ along an arc making one turn about $t=0$ in a positive direction. This shows that
$F^+(t)-F^-(t)$ is obtained as an analytic continuation of the function
$$F(t)-l_0^{*}(F(t)) = \frac{\int_{\delta(t)}\omega}{\int_{\delta(t)}\omega_0} -
\frac{\int_{l_0^*\delta(t)}\omega}{\int_{l_0^*\delta(t)}\omega_0}
=\frac{W_{\delta(t),l_0^*\delta(t)}(\omega,\omega_0)}{\int_{\delta(t)}\omega_0 \int_{l_0^*\delta(t)}\omega_0},
\quad t\in (-1/6,0)$$
along an arc contained in the upper complex half-plane. Here $ \omega=\mu_1 \omega_0+ \mu_2
\omega_{-1} + \mu_3 \omega_1 $, $l_0^*\delta(t),\delta(t)$ are cycles on the elliptic curve defined by $H_t=0$,
with two removed points $0^\pm$, and $l_0^*\delta(t) = \delta(t) - \gamma(t) +\alpha(t)$. Let $$i W(t),\quad t\in
(0,\infty)$$ be the real analytic function obtained as an analytic continuation of
$$
i W_{\delta(t),l_0^*\delta(t)}(\omega,\omega_0),\quad t\in (-1/6,0)
$$
along an arc contained in the upper  complex half-plane. As
$$
F^+(t)-F^-(t)= \frac{W(t)}{|\int_{\delta(t)}\omega_0|^2}
$$
then we shall show that $iW(t)$ has at most one zero on $(0,\infty)$. We use once again the argument principle.
We shall show that the analytic continuation of $W(t)$ in $\tilde{\mathcal{D}}=\C\setminus (-\infty,0)$ has at
most one zero counted with multiplicity. For this purpose we consider the complex domain
$$
\tilde{\mathcal{D}}_{Rr}= \tilde{\mathcal{D}}\cap \{ t: |t| > r, \;\; |t| < R \} .
$$
As before, denote by $W^\pm(t)$  the restriction on $(-\infty,0)$ of the analytic function obtained as an
analytic continuation of $W(t), t\in (0,\infty)$ along an arc contained in the upper (lower) complex half-plane
${\rm Im}\, \pm t
> 0$ respectively. Along the interval $(-1/6,0)$ we have
\begin{eqnarray*}
  W^+(t)-W^-(t) &=& W_{\delta(t),l_0^*\delta(t)}(\omega,\omega_0) - (l_0^*)^{-1}
W_{\delta(t),l_0^*\delta(t)}(\omega,\omega_0) \\
   &=& W_{\delta(t),l_0^*\delta(t)}(\omega,\omega_0) - W_{(l_0^*)^{-1}\delta(t),\delta(t)}(\omega,\omega_0) \\
   &=& W_{\delta(t),\delta(t)+\alpha(t)-\gamma(t)}(\omega,\omega_0) - W_{\delta(t)+\alpha(t)+\gamma(t),\delta(t)}(\omega,\omega_0) \\
   &=& 2  W_{\delta(t),\alpha(t)}(\omega,\omega_0)\\
   &=& -\frac{2}{\sqrt{2t}} \int_{\delta(t)}\omega_0 .
\end{eqnarray*}
Note that, according to section 5.1.1., the functions $W^\pm(t)$ are single-valued in a neighborhood of
$t=-1/6$. Therefore
$$
W^+(t)-W^-(t)= -\frac{2}{\sqrt{2t}} \int_{\delta(t)}\omega_0
$$
along the interval $(-\infty,0)$. We conclude that the imaginary part $i(W^+(t)-W^-(t)), t\in (-\infty,0)$ of
the analytic continuation of $iW(t),t\in (0,\infty)$ does not vanish.

It follows from the asymptotic expansion of $J_k$ near $t=0$ obtained in 5.3.1.(b) that
$$
W(t)\lesssim  \frac{1}{\sqrt{t}},\quad t\sim 0.
$$
Therefore along the small circle $\{ |t| = r \}$, the increase of the argument of $W(t)$ is at worst close to
$\pi$. Along the border of the big circle the decrease of the argument is close to $2\pi/6$, see $(5.8)$.
Summing up the above information we conclude that the total increase of the argument of $W(t)$ along the border
of $\tilde{\mathcal{D}}_{Rr}$ is close to (or smaller than) $4\pi-2\pi/6$ and hence $W(t)$ has at most one zero.
\emph{Therefore the maximal number of the zeros of $F(t)$ in $\mathcal{D}$ is two $($and this result is
exact$)$.  }

\vspace{2ex} \noindent
{\bf 5.4.2. Zeros of $F$ in the case (r11) and (5.2).}
In the same way as above we may study  $F(t)$ in the case (r11). This leads,
however, to a bound of the number of the zeros in $[-1/6,0)$ equal to three.
To improve the bound we consider first the function $\tilde{F}(t)=F(t)+\mu_4$,
where $\mu_4$ is a real constant.  We shall show that $\tilde{F}(t)$ has at
most three zeros in the complex domain $\mathcal{D}=\C\setminus[0,\infty)$.
For this we consider once again the domain $\mathcal{D}_{Rr}= \mathcal{D}\cap
\{t: r< |t| < R\}$ and evaluate the increase of the argument of $\tilde{F}$
along its border.  Along the boundary of the small circle $\{|t|=r\}$,
according to (5.11), the increase of the argument of $\tilde{F}(t)$ is close
to $3 \pi$, and along the boundary of the big circle $\{|t|=R\}$ this increase
is close to $0$, see (5.10).  Along the interval $(0,\infty)$ the imaginary
part of $F$ is  equal to
$$ \frac{W(t)}{|\int_{\delta(t)}\omega_0 |^2} $$
where, according to (5.5), for suitable constants $c_1,c_2$,
$$ W(t)=\frac{c_1t+c_2}{t^2}.$$
It follows that the imaginary part of $F(t)$
has at most one zero along $(0,\infty)$.  The argument principle implies that
$\tilde{F}(t)$ has at most three zeros in the complex domain $\mathcal{D}$.
Suppose now that for some $\mu_1,\mu_2,\mu_3$ the function $F(t)$ has exactly
three zeros in $\mathcal{D}$, and hence in the domain $\mathcal{D}_{Rr}$ for
$r$ sufficiently small and $R$ big enough.  According to (5.10) for $|t|$
sufficiently big we have $$ F(t) = c t^{-k/3} + o(|t^{-k/3}|) .  $$ It follows
that for sufficiently small $\mu_4$ a real zero of $\tilde{F}(t)$ bifurcates
from $\infty$ on the interval $(-\infty,0)$ and this zero is not contained in
the domain $\mathcal{D}_{Rr}$.  Therefore $\tilde{F}(t)$ will have at least
four zeros in the domain $\mathcal{D}$, in contradiction with the result proved
above.  \emph{Therefore the maximal number of the zeros of $F(t)$ in
$\mathcal{D}$ is two $($and this result is exact$)$.  }

Theorem 4 is proved. $\Box$

 \vspace{2ex} \noindent {\bf 5.5.  Proof of Theorem 3.}
 Denote the open period
annulus of the fixed reversible system (r18) or (r11) by $\Pi$.
Let $X_\lambda$, $\lambda\in \Lambda$ be the set of all quadratic plane
 vector fields, analytic with respect to $\lambda$ and such that $X_0$ coincides with
  (r18) (respectively, with (r11)). Theorem 3 can be reformulated as follows:

\vspace{1ex}
\emph{The cyclicity $Cycl(\Pi,X_\lambda)$ of the open period annulus $\Pi$ is equal to two } \\

\vspace{1ex}
\noindent
which means that $X_\lambda$ has {\it at most two} limit cycles
which tend to $\Pi$ as $\lambda$ tends to zero. It is shown in
\cite[Theorem 1]{gav1} that if the cyclicity  $Cycl(\Pi,X_\lambda)$ is finite, then there exists a germ of an
analytic curve $\lambda(\varepsilon)$, such that
$$
Cycl(\Pi,X_\lambda) = Cycl(\Pi,X_\lambda(\varepsilon)) .
$$
In other words, it is enough to study one-parameter deformations. Let $I(t)$ be the  first non-zero
Poincar\'{e}-Pontryagin (or generating) function associated to such a perturbation. The derivative $J(t)=I'(t)$  is
of the form $(5.12)$ or $(5.13)$ and we proved that this function has at most two zeros (Theorem 4), and hence
$I(t)$ has at most three zeros on $[-1/6,0)$, one of them being $t=-1/6$ . Therefore, by a standard argument,
the cyclicity of the open period annulus of the perturbed one-parameter quadratic system is at most two.

It remains to show that the cyclicity $Cycl(\Pi,X_\lambda)$ is finite. It is shown in \cite{gav1} that if
$Cycl(\Pi,X_\lambda)=\infty$, then there exists a Poincar\'{e}-Pontryagin function $I(t)$ (associated to some
one-parameter deformation) with infinite number of zeros. As $I(t)$ is an Abelian integral, then this is clearly
impossible. This completes the proof of Theorem 3. $\Box$

\vspace{2ex} \noindent {\bf Acknowledgment.}  This research has been partially
 supported by PAI  Rila program
through Grants  14749SM (France) and Rila 3/6-2006 (Bulgaria).

\vspace{2cm}

E-mail addresses:

 \vspace{1cm}
 sebagaut@wanadoo.fr

lubomir.gavrilov@math.univ-tlse.fr

iliya@math.bas.bg

\end{document}